\documentclass[12pt]{article}
\usepackage{amsmath,amsthm,amssymb,amscd,}
\usepackage{latexsym}
\newtheorem{thm}{Theorem}[section]

 \newtheorem{prop}[thm]{Proposition}
 \theoremstyle{definition}
 
 \theoremstyle{remark}

 \numberwithin{equation}{section}

\title
{Notes on branched coverings of Seifert manifolds
\thanks{Research supported by NSFC
(No. 11171025) and by Laboratory of Mathematics and Complex Systems, Ministry of Education.} }

\author{ Hong Huang
  }

\date{}
\begin{document}
\maketitle

\begin{abstract}  In a paper published in 2002, the author gave a criterion to determine whether there is a fiber-preserving branched covering with given degree between two given connected, closed, orientable Seifert manifolds with orientable bases. Here we supply some details of the proof of two claims in that paper.  We give  an explicit construction of fiber-preserving branched covering between two Seifert fibered solid tori  when their Seifert invariants satisfy certain relation, and we show the factorability  of fiber-preserving branched coverings between two closed Seifert manifolds.

\vspace{0.1in}

{\bf Key words}: Seifert manifolds, fiber-preserving branched coverings,  factorization of branched coverings

\vspace{0.1in} {\bf Mathematics Subject Classification 2010}: 57N10

\end{abstract}
\maketitle

 In a paper [2] published in 2002, the author gave a criterion to determine whether there is a fiber-preserving branched covering with given degree between two given connected,  closed, orientable  Seifert manifolds with orientable bases. There are two statements used in [2] without detailed proof: One is on the existence of  fiber-preserving branched covering between two Seifert fibered solid tori  when their Seifert invariants satisfy certain relation, another is on the factorability  of fiber-preserving branched coverings between two closed Seifert manifolds.
Her we'll supply  some details of the proof of these two statements.

 Let $V$ be a Seifert fibered solid torus. Let $(q, h)$ be a (cross section, fiber) basis of $H_1(\partial V)$. We will abuse notation to use $q, h, etc.$ to denote both a closed curve and the homology class that it represents. If a meridian $m$ of $V$ is homologous  to
 $\alpha q+\beta h$ in $H_1(\partial V)$, where $ \alpha>0, \beta$ are coprime integers, then we say that the   fibered solid torus $V$ has Seifert invariant $(\alpha, \beta)$ (w.r.t. $(q, h)$); see [5].

 The following result is contained in [2], but the details of the proof of sufficiency are omitted there.

\begin{prop}   Let $V_1, V_2$ be two Seifert fibered solid tori. Suppose the Seifert invariant of $V_2$ w.r.t. some (cross section, fiber) basis of $H_1(\partial V_2)$ is $(\alpha_2, \beta_2)$, and $k$ is a positive integer. Then there is a fiber-preserving $k$-fold branched covering $f:V_1\rightarrow V_2$ with fiber degree one and with upstairs and  downstairs branched set the central fiber of $V_1$ and  $ V_2$ respectively if and only if there is some (cross section, fiber) basis of $H_1(\partial V_1)$  such that the Seifert invariant $(\alpha_1, \beta_1)$ of $V_1$ w.r.t. this basis satisfies $\beta_1 / \alpha_1=k\beta_2 / \alpha_2$.
\end{prop}

\vspace{0.3cm}

\noindent {\bf Proof} \ \  First we show the necessity as in [2].  Let $f:V_1\rightarrow V_2$ be a fiber-preserving $k$-fold branched covering with fiber degree one and with upstairs and  downstairs branched set the central fiber of $V_1$ and  $ V_2$ respectively. (By upstairs branched set we mean the preimage of the downstairs branched set under $f$.) Let $(q_2, h_2)$ be the (cross section, fiber) basis of $H_1(\partial V_2)$ w.r.t. which the Seifert invariant of $V_2$ is $(\alpha_2, \beta_2)$. Let $q_1=f^{-1}(q_2)$. Since $\partial V_1$ is the union of fibers through all points of $q_1$  and the fiber degree of $f$ is one, $q_1$ must be connected as $\partial V_1$ is.
Then we see that $q_1$ is a cross section of the Seifert fibration of $\partial V_1$, and $f|_*(q_1)=kq_2$, where $f|$ is the restriction of $f$ to $\partial V_1$.

 Let the Seifert invariant of $V_1$ w.r.t. the (cross section, fiber) basis $(q_1, h_1)$ be $(\alpha_1, \beta_1)$. By assumption we have $f|_*(h_1)=h_2$. By definition of branched covering, there is a nonzero integer $s$ such that $f|_*(m_1)=sm_2$, where $m_1$ (resp. $m_2$) is a meridian of $V_1$ (resp. $V_2$).  So $f|_*(\alpha_1q_1+\beta_1h_1)=s(\alpha_2q_2+\beta_2h_2)$, and $\alpha_1kq_2+\beta_1h_2=s\alpha_2q_2+s\beta_2h_2$, which implies that $\alpha_1k=s\alpha_2$, and $\beta_1=s\beta_2$.  It follows that $\beta_1 / \alpha_1=k\beta_2 / \alpha_2$.

 Then we show the sufficiency. Suppose there is some (cross section, fiber) basis $(q_1,h_1)$ of $H_1(\partial V_1)$  such that the Seifert invariant $(\alpha_1, \beta_1)$ of $V_1$ w.r.t. this basis satisfies $\beta_1 / \alpha_1=k\beta_2 / \alpha_2$. If $\beta_2=0$, then $\alpha_2=1$, $\beta_1=0$, and $\alpha_1=1$. In this case, the result is clear. Below we assume that $\beta_2\neq 0$. From the relation $\beta_1 / \alpha_1=k\beta_2 / \alpha_2$, using that $\alpha_j$ and $\beta_j$ are coprime, we get that $\alpha_1|\alpha_2$, and $\beta_2|\beta_1$.

 Now it is clear that there exists a fiber-preserving $k$-fold branched covering $f: V_1\rightarrow V_2$  with fiber degree one and with upstairs and  downstairs branched set the central fiber of $V_1$ and  $ V_2$ respectively. The reason is as follows.  First we can construct a fiber preserving covering $f|:\partial V_1 \rightarrow \partial V_2$ such that $f|_*(q_1)=kq_2$, and $f|_*(h_1)=h_2$. We have, using the relation $\beta_1 / \alpha_1=k\beta_2 / \alpha_2$,
 \begin{equation*}
 f|_*(m_1)=f|_*(\alpha_1q_1+\beta_1h_1)=\alpha_1kq_2+\beta_1h_2=\frac{\beta_1}{\beta_2}(\alpha_2q_2+\beta_2h_2)=\frac{\beta_1}{\beta_2}m_2.
  \end{equation*}
  Then by shrinking $\partial V_1, \partial V_2$ and $f|$ respectively we see that $f|$ extends to a fiber-preserving branched covering $f:V_1\rightarrow V_2$ as desired. (Compare  the hint to Exercise 3 on p. 312 of [4].)

 Below we give a more explicit construction. Let $\alpha_j', \beta_j'$ be integers such that $\alpha_j \beta_j'-\alpha_j'\beta_j=1$, $j=1,2$; note that $\alpha_j'$ (resp.  $\beta_j'$) is determined by $\alpha_j$ and $\beta_j$ up to multiples of $\alpha_j$ (resp. $\beta_j$). Now we have
 $\alpha_1 \beta_1'-(\frac{\beta_1}{\beta_2} \alpha_1')\beta_2=1$, and $\alpha_1(\frac{\alpha_2}{\alpha_1} \beta_2')-\alpha_2'\beta_2=1$. (As noted above, both $\frac{\alpha_2}{\alpha_1}$ and $\frac{\beta_1}{\beta_2}$ are integers.) So $\alpha_1 (\beta_1'-\frac{\alpha_2}{\alpha_1} \beta_2')=(\frac{\beta_1}{\beta_2} \alpha_1'-\alpha_2')\beta_2$.  We see that    $\alpha_2'-\frac{\beta_1}{\beta_2}\alpha_1'$ is a multiple of $\alpha_1$, since   $\alpha_1$ and $\beta_2$ are coprime.

  For $j=1,2$, consider the  universal covering $p_j: D^2 \times \mathbb{R}\rightarrow V_j$ with deck transformation group generated by the element

 \vspace{0.3cm}

 $\begin{array}{l}
 \tau_j: D^2 \times \mathbb{R}\rightarrow D^2 \times \mathbb{R}, \\
 (re^{2\pi i \theta}, t) \mapsto  (re^{2\pi i(\theta-\frac{\alpha_j'}{\alpha_j})}, t+1).
 \end{array}$

 \vspace{0.3cm}

\noindent Here we think of $V_j$ as the quotient of $D^2 \times \mathbb{R}$ by the action of the corresponding deck transformation group; cf. [5].

First we construct a map $\tilde{f}: D^2 \times \mathbb{R}\rightarrow  D^2 \times \mathbb{R}$ as follows,
\begin{equation*}
\tilde{f}(re^{2\pi i \theta}, t)=(re^{2\pi i \frac{\beta_1}{\beta_2}\theta}, \frac{\alpha_2}{\alpha_1}t).
\end{equation*}

Note that

\vspace{0.3cm}

$\begin{array}{l}
\tilde{f}(\tau_1(re^{2\pi i\theta}, t)) \\
=\tilde{f}(re^{2\pi i(\theta-\frac{\alpha_1'}{\alpha_1})}, t+1) \\
=(re^{2\pi i\frac{\beta_1}{\beta_2}(\theta-\frac{\alpha_1'}{\alpha_1})}, \frac{\alpha_2}{\alpha_1}(t+1))\\
=(re^{2\pi i (\frac{\beta_1}{\beta_2}\theta-\frac{\alpha_2}{\alpha_1}\frac{\alpha_2'}{\alpha_2})}, \frac{\alpha_2}{\alpha_1}t+\frac{\alpha_2}{\alpha_1})\\
=\tau_2^{\alpha_2/\alpha_1}(re^{2\pi i \frac{\beta_1}{\beta_2}\theta}, \frac{\alpha_2}{\alpha_1}t) \\
=\tau_2^{\alpha_2/\alpha_1}(\tilde{f}(re^{2\pi i\theta}, t)).
\end{array}$

 \noindent Here in the third equality we have used the fact $\alpha_1|(\alpha_2'-\frac{\beta_1}{\beta_2}\alpha_1')$.

\noindent It follows that  $\tilde{f}$ descends to a map $f:V_1\rightarrow V_2$. One can easily check that $f: V_1\rightarrow V_2$ is the desired branched covering.
\hfill{$\Box$}

\vspace{0.3cm}

The following result  was stated in the proof of Theorem 2.2 in [2] without proof.

\begin{prop}  Let $f: M_1\rightarrow M_2$ be a fiber-preserving $d$-fold branched covering between two (connected) closed Seifert manifolds (both  orientable and with orientable base).  Then there exists a positive  integer $d_1$ with $d_1|d$ such that $f=f_2\circ f_1$, where $f_1:M_1\rightarrow M_1/\mathbb{Z}_{d_1}$ (here the $\mathbb{Z}_{d_1}$-action on $M_1$ is inside the $S^1$-action on $M_1$) is the natural projection, and $f_2: M_1/\mathbb{Z}_{d_1} \rightarrow M_2$ is a fiber-preserving  branched covering with fiber degree one.
\end{prop}
\noindent {\bf Proof} \ \ Let $\{O_1,\cdot\cdot\cdot,O_n\}$ be a nonempty collection of finitely many fibers
in $M_2$, including all singular fibers of the Seifert fibration of $M_2$ and downstairs branch set of $f$, and such that all the singular fibers of the Seifert fibration of $M_1$ are contained in $f^{-1}(O_1\cup \cdot\cdot\cdot \cup O_n)$.  Let $V_1, \cdot\cdot\cdot,  V_n$ be pairwise disjoint Seifert fibered solid tori in $M_2$ with central fibers $O_1,\cdot\cdot\cdot,O_n$ respectively, and $N_2 = M_2\setminus \text {int} (V_1 \cup \cdot\cdot\cdot \cup V_n)$.  Let $N_1=f^{-1}(N_2)$.  Then both $N_1$ and $N_2$ are trivial $S^1$-bundles (over compact, connected, orientable surfaces with boundary), and $f|N_1: N_1\rightarrow N_2$ is a fiber-preserving $d$-fold unbranched covering.  When restricted to a fiber of $N_1$, $f$ is a $d_1$-fold covering of a fiber in $N_2$ for some positive integer $d_1$ with $d_1|d$. Now we see that $f|N_1=f
_2'\circ f_1'$, where $f_1':N_1\rightarrow N_1/\mathbb{Z}_{d_1}$ (here the $\mathbb{Z}_{d_1}$-action on $N_1$ is inside the $S^1$-action on $N_1$ coming from that on $ M_1$) is the natural projection, and $f_2': N_1/\mathbb{Z}_{d_1} \rightarrow N_2$ is a fiber-preserving  unbranched covering with fiber degree one.  (By the way, one can easily check  that  $f_1'$ (resp.  $f_2'$) sends a meridian in any boundary torus of $N_1$ (resp. $N_1/\mathbb{Z}_{d_1}$) to a multiple of a meridian in the image torus.) Now we let $V_1, \cdot\cdot\cdot,  V_n$ be shrunk to the central fibers $O_1,\cdot\cdot\cdot,O_n$ (at the same time), and let $f_1$ (resp. $f_2$) be the limit of the mappings $f_1'$ (resp. $f_2'$) thus gotten in the process.  Then we are done.   (Compare [1].)   \hfill{$\Box$}

\vspace{0.3cm}

\noindent {\bf Remark}. \ \  See Theorem 3.3.1 in [3] for a result similar to Proposition 0.2 above. Our proof here is slightly  different, and is along the lines of proof of Lemma 2.1 in [2].

\vspace*{0.4mm}
 School of Mathematical Sciences, Beijing Normal University,

 Laboratory of Mathematics and Complex Systems, Ministry of Education,

 Beijing 100875, P.R. China

 E-mail address: hhuang@bnu.edu.cn

\end{document}